\documentclass[11pt,a4paper]{article}
\usepackage{ifthen,latexsym,amssymb,amsmath,bbm,fixmath,amsthm}
\usepackage[shortlabels]{enumitem} 

\newcommand{\I}[1]{{\mathbbm #1}}

\newcommand{\floor}[1]{\lfloor #1\rfloor}
\newcommand{\me}{{\mathrm e}}
\renewcommand{\mid}{:}

\renewcommand{\ge}{\geqslant}
\renewcommand{\le}{\leqslant}

\newif\ifnotesw\noteswtrue

\newcommand{\hide}[1]{}

\newcommand{\beq}[1]{\begin{equation}\label{#1}}
\newcommand{\eeq}{\end{equation}}

\newtheorem{theorem}{Theorem}[section]
\newtheorem{lemma}[theorem]{Lemma}

\newtheorem{corollary}[theorem]{Corollary}
\theoremstyle{remark}
\newtheorem{remark}[theorem]{Remark}

\newcommand{\bpf}[1][Proof.]{\smallskip\noindent{\it #1} }
\newcommand{\epf}{\qed \medskip}

\author{Oleg Pikhurko\footnote{Supported 
by ERC Advanced Grant 101020255.}\\
Mathematics Institute and DIMAP\\
University of Warwick\\
Coventry CV4 7AL, UK
}

\begin{document}

\title{Constructions of Tur\'an systems that are tight up to a multiplicative constant}
\maketitle

\begin{abstract}
For positive integers $n\ge s> r$, the \emph{Tur\'an function} $T(n,s,r)$ is the smallest size of an $r$-graph with $n$ vertices such that every set of $s$ vertices contains at least one edge. Also, define the \emph{Tur\'an density} $t(s,r)$ as the limit of $T(n,s,r)/ {n\choose r}$ as $n\to\infty$. The question of estimating these parameters received a lot of attention after it was first raised by Tur\'an in~1941. A trivial lower bound is $t(s,r)\ge 1/{s\choose s-r}$. In the early 1990s, de Caen conjectured that $r\cdot t(r+1,r)\to\infty$  as $r\to\infty$ and offered 500 Canadian dollars for resolving this question.
 
We disprove this conjecture by showing more strongly that for every integer $R\ge1$ there is $\mu_R$ (in fact, $\mu_R$ can be taken to grow as $(1+o(1))\, R\ln R$) such that 
 $t(r+R,r)\le (\mu_R+o(1))/ {r+R\choose R}$ as $r\to\infty$, that is, the trivial lower bound is tight for every $R$ up to a multiplicative constant $\mu_R$. 
\end{abstract}


\newcommand{\extend}[4]{#1\otimes_{#3}{#2}_{*}^{#4}}
\newcommand{\Ball}{\mathcal{B}}

\section{Introduction}

An \emph{$r$-graph $G$} on a vertex set $V$ is a collection of $r$-subsets of $V$, called \emph{edges}. If the vertex set $V$ is understood, then we may identify $G$ with its edge set, that is, view $G$ as a subset of  ${V\choose r}:=\{X\subseteq V\mid |X|=r\}$. 

For positive integers $n\ge s> r$, a \emph{Tur\'an $(n,s,r)$-system} is an $r$-graph $G$ on an $n$-set $V$ such that every $s$-subset $X\subseteq V$ is \emph{covered} by $G$, that is, there is $Y\in G$ with $Y\subseteq X$.  
Let the \emph{Tur\'an number} $T(n,s,r)$ be the smallest size of a Tur\'an $(n,s,r)$-system. If we pass to the complements then ${n\choose r}-T(n,s,r)$ is $\mathrm{ex}(n,K_s^r)$, the maximum size of an $n$-vertex $r$-graph without $K_s^r$, the complete $r$-graph on $s$ vertices. This is a key instance of the classical extremal \emph{Tur\'an problem} that goes back to Tur\'an~\cite{Turan41}.  For the purposes of this paper, it is more convenient to work with the function $T$ rather than with its complementary version $\mathrm{ex}$. So we define the \emph{Tur\'an density} 
 \beq{eq:t}
 t(s,r):=\lim_{n\to\infty} \frac{T(n,s,r)}{{n\choose r}},
 \eeq
 to be the asymptotically smallest edge density of a Tur\'an $(n,s,r)$-system as $n\to\infty$.
The limit in the right-hand size of~\eqref{eq:t} exists, since easy double-counting shows that $T(n,s,r)/{n\choose r}$ is non-decreasing in $n$.
 
 We refer the reader to de Caen~\cite{Decaen94} and Sidorenko~\cite{Sidorenko95} for  surveys of results and open problems on minimum Tur\'an systems, and to F\"uredi~\cite{Furedi91} and Keevash~\cite{Keevash11} for surveys of the $\mathrm{ex}$-function for general hypergraphs. 
 Also, the table on Page 651 in Ruszink\'o~\cite{Ruszinko07} lists some connections of Tur\'an systems to various areas of combinatorial designs.


In the trivial case $r=1$, it holds that $T(n,s,1)=n-s+1$ for any $n\ge s> 1$. The case $r=2$ was resolved in the fundamental paper of Tur\'an~\cite{Turan41} (with the special case when $(s,r)=(3,2)$ previously done by Mantel~\cite{Mantel07}). In particular, it holds that $t(s,2)=\frac1{s-1}$ with the upper bound coming from the disjoint union of $s-1$ almost equal cliques. 
The problem of determining $T(n,s,r)$ for $r\ge 3$ was raised already in the above-mentioned paper of Tur\'an~\cite{Turan41} from 1941. 
Erd\H os~\cite[Section~III.1]{Erdos81} offered \$500 for determining $t(s,r)$ for at least one pair $(s,r)$ with $s>r\ge 3$. This reward is still unclaimed despite decades of active research. It was conjectured by Tur\'an and other researchers (see e.g.\ \cite[Page 348]{Decaen83}) that $t(s,3)=4/({s-1})^2$ for each $r\ge 4$. Various constructions attaining this upper bound can be found in Sidorenko's survey~\cite[Section 7]{Sidorenko95}. 
In the first open case $s=4$, 
the computer-generated lower bound  $t(4,3)\ge 0.438...$ of Razborov~\cite{Razborov10} (improving on the earlier bounds in~\cite{KatonaNemetzSimonovits64,Decaen88,ChungLu99}) comes rather close to the conjectured value $4/9$. As stated by Sidorenko~\cite[Section~8]{Sidorenko95} (and this seems to be still true), the only pair $(s,r)$ with $s>r\ge 4$ for which there is a plausible conjecture is $(5,4)$, where a construction of Giraud~\cite{Giraud90} gives $t(5,4)\le \frac5{16}=0.325$. (The best known lower bound is $t(5,4)\ge \frac{627}{2380}=0.2634...$ by
Markstr\"om~\cite{Markstrom09}.)

Also, a lot of attention was paid to estimating $t(s,r)$ as $r\to\infty$. In the first interesting case $s=r+1$, the trivial lower bound $t(r+1,r)\ge \frac1{r+1}$ was improved to $1/r$ independently by de Caen~\cite{Decaen83ac}, Sidorenko~\cite{Sidorenko82}, and Tazawa and Shirakura~\cite{TazawaShirakura83}. Some further improvements (of order at most $O(1/r^2)$) for a growing sequence of $r$  were made by Giraud (unpublished, see~\cite[Page 362]{ChungLu99}), Chung and Lu~\cite{ChungLu99}, and by Lu and Zhao~\cite{LuZhao09}. 

In terms of the previously known upper bounds on $t(r+1,r)$ as $r\to\infty$ there was a sequence of better and better bounds: $O(1/\sqrt r)$ by Sidorenko~\cite{Sidorenko81}, $\frac{1+2\ln r}{r}$ by Kim and Roush~\cite{KimRoush83}, $\frac{\ln r+O(1)}{r}$ by Frankl and R\"odl~\cite{FranklRodl85gc}, and  $(1+o(1))\frac{\ln r}{2r}$ by Sidorenko~\cite{Sidorenko97}.

De Caen~\cite[Page 190]{Decaen94} conjectured that $r\cdot t(r+1,r)\to\infty$ as $r\to\infty$ and offered 500 Canadian dollars for proving or disproving this. The question which of the bounds $\Omega(\frac 1r)\le t(r+1,r)\le O(\frac{\ln r}{r})$ is closer to the truth  was asked earlier by Kim and Roush~\cite[Page 243]{KimRoush83}.
Also, Frankl and R\"odl~\cite[Page~216]{FranklRodl85gc} wrote that it is ``conceivable" that $t(r+1,r)=O(1/r)$.  
Here we disprove de Caen's conjecture 
(and thus confirm the intuition of Frankl and R\"odl), 
with the following explicit constants.

\begin{theorem}\label{th:main} For all integers $n> r\ge 1$, it holds that $T(n,r+1,r)\le \frac{6.239}{r+1}\,{n\choose r}$. 

Also, there is $r_0$ such that, for all integers $n>r\ge r_0$, it holds that $T(n,r+1,r)\le \frac{4.911}{r+1}\,{n\choose r}$.\end{theorem}

While the constant in the first part is worse than in the second part, we include both proofs as it may be useful to have a simple explicit bound valid for every pair $(n,r)$. 

In the general case, the trivial lower bound is $T(n,s,r)\ge {n\choose r}/{s\choose s-r}$: indeed, if $G\subseteq {[n]\choose r}$ has fewer edges than the stated bound then the expected number of edges inside a random $s$-set is ${s\choose r}\cdot |G|/{n\choose r}<1$ so some $s$-set is not covered at all. Thus $t(s,r)\ge 1/{s\choose s-r}$. This was improved to 
\beq{eq:Decaen83}
 t(s,r)\ge \frac{1}{{s-1\choose s-r}},
\eeq
 by de Caen~\cite{Decaen83}, and this is still  the best known general lower bound. In terms of upper bounds, Frankl and R\"odl~\cite{FranklRodl85gc} proved that, for any integer $R\ge 1$, we have 
 \beq{eq:FR}
 t(r+R,r)\le (1+o(1)) \frac{R(R+4)\ln r}{{r+R\choose R}},\quad\mbox{as $r\to\infty$}.
 \eeq

We can also remove the factor $\ln r$ in the above result of Frankl and R\"odl:

\begin{theorem}\label{th:New} 
For every integer $R\ge 1$ it holds that
 \beq{eq:new}
 t(r+R,r)\le (\mu+o(1))\,\frac{1}{{r+R\choose R}},\quad\mbox{as $r\to \infty$},
 \eeq
 where $\mu:=(c_0+1)^{R+1}/c_0^{R}$ with $c_0=c_0(R)$ being the largest real root of the equation $\me^{c}=(c+1)^{R+1}$.
  \end{theorem}

While for a given integer $R$ it is possible to numerically approximate the above constant $\mu$ (in particular, to see that $\mu<4.911$ for $R=1$), it is also interesting to see how $\mu$ grows with~$R$. This is done in the following corollary to Theorem~\ref{th:New}:
   
\begin{corollary}\label{cr:General}   
 There is $R_0$ such that for every integer $R\ge R_0$, it holds that
 \beq{eq:cr:General}
 \limsup_{r\to\infty}\, t(r+R,r)\cdot {r+R\choose R}\le R\ln R + 3R\ln\ln R.
\eeq
\end{corollary}

Let us discuss the known upper bounds in the case when $r\to\infty$ while $R=R(r)$ is a function of $r$ that also tends to the infinity. By analysing the construction of Frankl and R\"odl~\cite{FranklRodl85gc}, Sidoren\-ko~\cite[Theorem~2]{Sidorenko97} proved that $\mu(r+R,r)\le (1+o(1)) R \ln {r+R\choose r}$ 
provided $R\ge r/\log_2 r$. Liu and Pikhurko~\cite{LiuPikhurko25x} observed that the restriction on $R$ in Sidorenko's bound can be weakened to just $R\to\infty$.
Also, Liu and Pikhurko~\cite{LiuPikhurko25x} calculated that the recursive construction presented in this paper yields $\mu(r+R,r)\le (1+o(1)) R\ln R$ for any function $R=o(\sqrt{r})$.

\section{Proofs}

Our construction is motivated by the recursive constructions of covering codes in~\cite{CooperEllisKahng02,KrivelevichSudakovVu03,LenzRashtchianSiegelYaakobi21}, of which Theorems~6 and 9 in Lenz, Rashtchian, Siegel and Yaakobi~\cite{LenzRashtchianSiegelYaakobi21} are probably closest  to the presented results. This connection was previously used by Verbitsky and Zhukovskii (personal communication) to prove new results on insertion covering codes using some methods developed for Tur\'an systems; this project later developed into a joint paper~\cite{PikhurkoVerbitskyZhukovskii}. Here, we exploit this connection by using the recursion in~\cite{LenzRashtchianSiegelYaakobi21} as guiding intuition for our construction.

It will be convenient to extend the definition of $T(n,s,r)$ to allow all triples $(n,s,r)$ of integers with $s>r\ge 0$ and $n\ge 0$. 
We agree (and this formally matches the general definition) that  $T(n,s,r)=0$ for $n<s$ (in particular, for $n\le r$), while $T(n,s,0)=1$ 
for $n\ge s\ge 1$. These degenerate cases  will be used in our inductive proofs.

We need some definitions first.
For integers $0\le m\le n$, we denote $[n]:=\{1,2,\dots,n\}$ and $[m,n]:=\{m,m+1,\dots,n\}$.
 For an integer $m\ge 0$ and an $\ell$-graph $H\subseteq {[n]\choose \ell}$, let $\extend{H}{K}{n}{m}$ denote the $(\ell+m)$-graph on $[n]$ that consists of those $X\in {[n]\choose \ell+m}$ for which there is $Y\in H$ with $Y$ being the initial $\ell$-segment of $X$, that is, if we order the elements of $X$ as $x_1<\dots<x_{\ell+m}$ under the natural ordering of $[n]$ then $\{x_1,\dots,x_\ell\}\in H$. 
Informally speaking, $\extend{H}{K}{n}{m}$ is obtained from $H$ by extending its edges into $(\ell+m)$-subsets of $[n]$  in all possible ways to the right. 
Also, we define
$$
\Ball(H):=\left\{B\in{[n]\choose \ell+1}\mid {B\choose \ell}\cap H\not=\emptyset\right\}, 
$$
to consist of all $(\ell+1)$-subsets of $[n]$ covered by the $\ell$-graph~$H$.

The first part of Theorem~\ref{th:main} will be derived from the following lemma.

\begin{lemma}\label{lm:main}
Let reals $\beta\in (0,1)$ and $c,\mu>1$ be fixed such that 
$\floor{\beta \mu}\ge c$
and
 \beq{eq:mu1}
 \frac{c}{\beta\mu-1}+\frac{\me^{-c}}{1-\beta}\le 1.
\eeq
 Then, for all integers $n,r\ge 0$, there is a Tur\'an $(n,r+1,r)$-system $G_n^r$ with $|G_n^r|\le \frac{\mu}{r+1}{n\choose r}$. 
\end{lemma}

\bpf We construct $G_n^r$ using induction on $r$ and then on~$n$. Let $r_0:=\floor{\mu}-1$.

For $r\le r_0$ and any $n\ge 0$, we can let $G_n^r:={[n]\choose r}$ be the complete $r$-graph on $[n]$. Note that $\frac{\mu}{r+1}\ge 1$ so the desired upper bound $|G_n^r|\le \frac{\mu}{r+1}{n\choose r}$ trivially holds.

Let $r>r_0$. Given $r$, we construct $G_n^r$ inductively on $n$. For $n\in [0,r]$, we let $G_{n}^r:=\emptyset$ be the empty $r$-graph on~$[n]$, which trivially satisfies the lemma. 

Let $n\ge r+1$. 
Define $k:=\floor{\beta(r+1)}$. Note that $k\ge 1$ since $r+1\ge r_0+2> \mu$ while $\beta \mu\ge c> 1$ by our assumptions. Also, $k\le r$ since $\beta<1$. We will consider a random $r$-graph $G_n^r$ (which will be a Tur\'an $(n,r+1,r)$-system deterministically) and fix an outcome whose size is at most the expected value.

Let $S\subseteq {[n]\choose k-1}$ be a \emph{$\frac{c}{k}$-binomial} random $(k-1)$-graph on $[n]$, that is, we include each $(k-1)$-subset of $[n]$ into $S$ with  probability ${c}/{k}$, with all choices being mutually independent. (Note that $k\ge \floor{ \beta(r_0+2)}\ge \floor{\beta\mu}$ which is at least $c$ by one of the assumptions, so $c/k\le 1$.) 
The expected size of the $r$-graph $S^*:=\extend{S}{K}{n}{r-k+1}$  is exactly $\frac ck\, {n\choose r}$ because every $r$-set $Y\in \binom{[n]}{r}$ is included into it with probability $c/k$: indeed, $Y\in S^*$ if and only if the initial $(k-1)$-segment of $Y$ is in $S$, which happens with probability~$c/k$. 

Let $T:={[n]\choose k}\setminus \Ball(S)$, that is, $T$ consists of those $k$-subsets $Y$ of $[n]$ such that no $(k-1)$-subset of $Y$ belongs to~$S$. Note that
every $Y\in {[n]\choose k}$ is included into $T$ with probability exactly $\left(1-\frac ck\right)^k$: each of its $k$ subsets of size $k-1$ has to be omitted from~$S$. Let $T^*:=\extend{T}{G}{n}{r-k}$ be the $r$-graph on $[n]$ obtained as follows: for every edge $Y\in T$, let $y:=\max Y$, take a copy $G_Y$ of $G_{n-y}^{r-k}$ on $[y+1,n]$ and add to $T^*$ all sets $Y\cup Z$ with $Z\in G_Y$. (Note that if $y\ge n-r+k$ then no edges are added to $T^*$ for this $Y$.) By the inductive assumptions (since $k\ge 1$) and by $1-x\le \me^{-x}$, the expected size of $T^*$ can be upper bounded as follows, where $y$ plays the role of $\max Y$ for $Y\in T$:
 \begin{eqnarray*}
 \I E |T^*|&=& \sum_{y=k}^{n}  \left(1-\frac ck\right)^k{y-1\choose k-1}\cdot |G_{n-y}^{r-k}|
\\ &\le&  \sum_{y=k}^{n}  \me^{-c}{y-1\choose k-1} \cdot  \frac{\mu}{r-k+1}\, {n-y\choose r-k} 
\\ &= &  \frac{\me^{-c}\,\mu}{r-k+1}\, {n\choose r}
.
\end{eqnarray*}
 Note that we may have $k=r$ or $y=n$ in the above expressions; this is why we found  is convenient to allow any $n,r\ge 0$ when defining Tur\'an $(n,s,r)$-systems.
 
 Fix $S$ such that $|S^*\cup T^*|$ is at most its expected value, and let $G_n^r:=S^*\cup T^*$. We have by above that
 \beq{eq:G}
 |G_n^r|\le \I E|S^*\cup T^*|\le  \I E|S^*|+ \I E|T^*|\le
 \left(\frac{c}{k} + \frac{\me^{-c}\,\mu}{r-k+1}\right) {n\choose r}.
 \eeq
 Since  $r\ge r_0+1\ge \floor{\mu}$ and thus $r+1\ge \mu$, we can lower bound $k$ as
 $$
  k\ge {\beta(r+1)-1}= {(r+1)\left(\beta-\frac1{r+1}\right)}\ge  {(r+1)\left(\beta-\frac1{\mu}\right)} =\frac{(r+1)(\beta\mu-1)}{\mu}. 
 $$
  Using this and the trivial bound $r-k+1\ge (r+1)(1-\beta)$, we obtain from~\eqref{eq:G} that
 \[
  |G_n^r|\le  \left(\frac{c\mu}{(r+1)(\beta\mu-1)} + \frac{\me^{-c}\,\mu}{(r+1)(1-\beta)}\right) {n\choose r},
\]
 which is at most the claimed bound $\frac{\mu}{r+1}\, {n\choose r}$ by~\eqref{eq:mu1}.

Let us check that, regardless of the choice of $S$, the obtained $r$-graph $G_n^r$ is a Tur\'an $(n,r+1,r)$-system. Take any $(r+1)$-subset $X$ of $[n]$. Let its elements be $x_1<\dots<x_{r+1}$ and let $Y:=\{x_1,\dots,x_k\}$. If $Y$ is in $T$, then the Tur\'an $(n-x_k,r-k+1,r-k)$-system $G_Y$ on $[x_{k}+1,n]$ contains an edge $Z$ which is a subset of $\{x_{k+1},\dots,x_{r+1}\}\in{[x_{k}+1,n]\choose r-k+1}$; thus $G_n^r$ contains $Y\cup Z$ which is a subset of $X$, as desired. So suppose that $Y$ is not in $T$. By definition, this means that there is $i\in [k]$ such that $Z:=Y\setminus \{x_i\}$ is in~$S$. But then $X \setminus \{x_i\}$ has $Z$ as its initial $(k-1)$-segment and thus belongs to $S^*\subseteq G_n^r$, as desired. This finishes the proof of the lemma.\epf

It is easy to find a triple $(\beta,c,\mu)$ satisfying Lemma~\ref{lm:main}: for example, take $\beta=1/2$, $c=2$ and sufficiently large $\mu$. The constant in the first part of Theorem~\ref{th:main} is the optimal $\mu$ coming from Lemma~\ref{lm:main}, rounded up in the third decimal digit.

\bpf[Proof of the first part of Theorem~\ref{th:main}.] The assignment $\beta:=0.784$, $c:=2.89$ and $\mu:=6.239$ can be checked to satisfy Lemma~\ref{lm:main}.%
\hide{\begin{verb}
p = m /. Solve[c/(b*m - 1) + E^(-c)/(1 - b) == 1, m][[1]]
p /. b -> 0.784309 /. c -> 2.893122
p /. b -> 0.784 /. c -> 2.89
p /. b -> 0.78 /. c -> 2.89

p = m /. Solve[c/b + m*E^(-c)/(1 - b) == m, m][[1]]
p /. b -> 0.71533 /. c -> 2.51286
p /. b -> 0.715 /. c -> 2.51
p /. b -> 0.71 /. c -> 2.51
\end{verb}}
\epf

\begin{remark} The best constant $\mu=6.2387...$ for Lemma~\ref{lm:main} comes from some small values of~$r$. It can be improved in various ways, even just by using $T(n,3,2)\le \frac12{n\choose 2}$ in the base case of the induction in the proof.\end{remark}

Next we turn to general $R$ (with the result also including the case $R=1$, which will be used to derive the second part of Theorem~\ref{th:main}).

\begin{lemma}\label{lm:General}
If an integer $R\ge 1$, and reals $\beta\in (0,1)$ and $c,\mu>0$  satisfy $\me^{-c}<(1-\beta)^R$ and
 \beq{eq:General}
  \frac{c}{\beta^R}+\frac{\me^{-c}\mu}{(1-\beta)^R}\le \mu,
  \eeq
 then there is a constant $D$ such that, for all integers $n,r\ge 0$, there is a Tur\'an $(n,r+R,r)$-system $H_n^r$ with  
  \beq{eq:GeneralAim}
   |H_n^r|\le \left(\mu + \frac D{\ln(r+3)}\right) {r+R\choose R}^{-1}{n\choose r}.
   \eeq
 \end{lemma}
 
\bpf Given $R$, $\beta$, $c$ and $\mu$, fix constants $C$, $r_0$ and $D$ in this order, with each being sufficiently large depending on the previous constants (and with $r_0$ being an integer). We construct $H_n^r\subseteq {[n]\choose r}$ by induction on $r$ and, for each $r$, by induction on $n$. For $r\in [0,r_0]$ and any $n$, we can take the complete $r$-graph on $[n]$ for $H_n^r$. Note that~\eqref{eq:GeneralAim} holds since we can assume that ${r+R\choose R}\le D/\ln (r+3)$ for every such $r$. (We use that $\ln(r+3)>0$ for every $r\ge 0$, which was the only reason why $3$ was added to the argument of $\ln$.)

So let $r>r_0$. Define $k:=\floor{\beta r}$. We have that $k\ge \floor{\beta r_0}$ is at least $R$ because $r_0$ is sufficiently large. 
For $n\in [0,r]$, we let $H_n^r$ be the empty $r$-graph on $[n]$, which trivially satisfies the lemma. So let $n\ge r+1$.

Let $S$ be a random subset of ${[n]\choose k-R}$ where each $(k-R)$-subset of $[n]$ is included with probability $c/{k\choose R}$ independently of the others. Note that $c/{k\choose R}\le 1$
since $r\ge r_0$ is sufficiently large depending on $\beta$, $c$ and $R$. Let $S^*:=\extend{S}{K}{n}{r-k+R}$. Recall that this is the $r$-graph obtained by extending the edges of $S$ to the right into all possible $r$-subsets of~$[n]$. Also, let $T:={[n]\choose k}\setminus\Ball_R(S)$, where 
$$
\Ball_R(S):=\left\{X\in {[n]\choose k}\mid {X\choose k-R}\cap S\not=\emptyset\right\}
$$ 
consists of all $k$-subsets of $[n]$ covered by~$S$. Let $T^*:=\extend{T}{H}{n}{r-k}$ be the $r$-graph on $[n]$ obtained as follows.
For every edge $Y$ in $T$, let $y:=\max Y$, take a copy $H_Y$ of
$H_{n-y}^{r-k}$ on $[y+1,n]$  and, for every $Z\in H_Y$, add $Y\cup Z$ to~$T^*$. 
Take $S$ such that the size of $H_n^r:=S^*\cup T^*$ is at most its expected value.

Let us show that  $H^r_n$ is a Tur\'an $(n,r+R,r)$-system, regardless of the choice of~$S$. Take any $(r+R)$-subset $X\subseteq [n]$ with elements $x_1<\dots<x_{r+R}$. Let $Y:=\{x_1,\dots,x_k\}$. If $Y$ is in $T$ then some edge $Z$ in the Tur\'an $(n-x_k,r-k+R,r-k)$-system $H_Y$ satisfies $Z\subseteq X\setminus Y$ and thus $Y\cup Z\in T^*$ is a subset of~$X$. Otherwise there is an $R$-subset $Z\subseteq Y$ such that $Y\setminus Z\in S$; then $X\setminus Z$ is in $S^*$ and is a subset of~$X$. Thus every $(r+R)$-subset of $[n]$ is covered by $H_n^r$, as desired.

It remains to show that $H_n^r$ satisfies~\eqref{eq:GeneralAim}. Similarly as in Lemma~\ref{lm:main}, we have by induction that
 \begin{eqnarray*}
 |H_n^r|&\le& \I E|S^*|+\I E|T^*|
\\&\le& \frac{c}{{k\choose R}} {n\choose r}
+\sum_{y=k}^{n}  
\left(1-c\,{k\choose R}^{-1}\right)^{{k\choose R}}
{y-1\choose k-1} \cdot |H_{n-y}^{r-k}|
\\&\le& \left(\frac{c}{{k\choose R}}+  \frac{\me^{-c}}{{r-k+R\choose R}}\left(\mu+\frac{D}{\ln(r-k+3)}\right)\right){n\choose r}.
 \end{eqnarray*}
 Thus we have 
  \begin{eqnarray*}
  |H_n^r|\,\frac{ {r+R\choose R}}{{n\choose r}} 
   &\le &\frac{c}{\beta^R}+\frac{\me^{-c}\,\mu}{(1-\beta)^R}+\frac Cr+\left(\left(\frac{\me^{-c}}{(1-\beta)^R}+\frac Cr\right) \frac{D}{\ln (r-k+3)}\right)
    \\&\le & \mu+ \frac{D}{\ln (r+3)} + \left(\left(\frac{\me^{-c}}{(1-\beta)^R}-1+\frac1C\right) \frac{D}{\ln(r-k+3)}\right)
  \ \le\ \mu+ \frac{D}{\ln (r+3)},
  \end{eqnarray*}
  as desired. Here, the first inequality uses the fact that the ratio ${r+R\choose R}/{k\choose R}$ (resp.\ ${r+R\choose R}/{r-k+R\choose R}$) deviates from the ``main'' term $\beta^{-R}$ (resp.\ $(1-\beta)^{-R})$ by $O(1/r)$, and the resulting error can be absorbed by $C/r$. The second inequality uses~\eqref{eq:General} and absorbs all error terms by the (much larger) term $D/(C\ln(r-k+3))$. The last inequality uses the assumption that $\me^{-c}<(1-\beta)^R$.
  
   This finishes the proof of the lemma.\epf 

\bpf[Proof of Theorem~\ref{th:New}.] It is enough to show that if we take $c:=c_0$ as in the theorem  (that is, the largest root of $\me^{c}=(c+1)^{R+1}$), $\beta:=\frac{c}{1+c}$ and $\mu:=(c+1)^{R+1}/c^{R}$ then all assumptions of Lemma~\ref{lm:General} are satisfied. Note that $\me^{c}-(c+1)^{R+1}$ is negative (resp.\ positive) for sufficiently small (resp.\ large) $c>0$ so $c_0$ is well-defined and positive.

In fact, this assignment was obtained as follows: first, we solved~\eqref{eq:General} as equality for $\mu$, obtaining that $\mu=h(\beta,c)$, where
 \beq{eq:h}
 h(\beta,c):= \frac{c}{\beta^{R}-\me^{-c}(\frac{\beta}{1-\beta})^{R}},
 \eeq
 and then took  a  point at which the partial derivatives $\frac{\partial h}{\partial c}$ and $\frac{\partial h}{\partial \beta}$ vanish.%
 \hide{\begin{verbatim}
 h = (c/b^r)/(1 - E^(-c)/(1 - b)^r)
 
 hb = Simplify[D[h, b]]
 
 hb0 = 1 - (1 - b)^(1 + r) E^c
 
 hb = Simplify[D[h, c]]
 
 hc0 = -1 - c + (1 - b)^r E^c
 
 Simplify[hb0 + (1 - b) hc0]
 
 h00c = Simplify[h /. b -> c/(1 + c)]
 
 Simplify[hc0 /. b -> c/(1 + c)]
 \end{verbatim}}
 (It seems that this choice of $(\beta,c)$ minimises $\mu$ over the feasible region; however, this is not needed in our proof.)

Let us check that all assumptions of Lemma~\ref{lm:General} are satisfied. For $h$ as in~\eqref{eq:h} (and $c=c_0$), we have
 $$
 h\left(\frac{c}{1+c},c\right) = \frac{c}{(\frac{c}{1+c})^R-\me^{-c} c^R}=\frac{c}{(\frac{c}{1+c})^R- (\frac{1}{c+1})^{R+1}c^R}=\frac{(c+1)^{R+1}}{c^R}=\mu.
 $$
 This calculation also verifies the other constraint $(1-\beta)^R>\me^{-c}$, which is equivalent to the positivity of the denominator of~$h$.%
\epf

\bpf[Proof of the second part of Theorem~\ref{th:main}.] The constant $\mu$ given by Theorem~\ref{th:New} for $R=1$ can be seen to be $4.9108...$\,
, which is less than the constant in the stated upper bound. Alternatively, it is enough just to give some feasible value of $(\beta,c)$ such that the function $h(\beta,c)$ defined in~\eqref{eq:h} is less than $4.911$; one can check that a pair $(0.715,2.51)$ satisfies this.\epf

\bpf[Proof of Corollary~\ref{cr:General}.] Let $R$  be sufficiently large and let $c_0$ be as defined in Theorem~\ref{th:New}, that is, $c_0$ is the largest root of the equation $\me^c=(c+1)^{R+1}$. 

Let us show that
  \beq{eq:c0}
  R\ln R+R\ln\ln R< c_0< R\ln R+2R\ln\ln R.
  \eeq

For the upper bound, we have to show that, for any $c\ge R\ln R+2R\ln\ln R$, it holds that $\me^c>(c+1)^{R+1}$, or by taking the logarithms that $c>(R+1)\ln(c+1)$. If, say, $c+1\le \me R\ln R$ then 
  \begin{eqnarray*}
   c-(R+1)\ln(c+1) &\ge& R\ln R+2R\ln\ln R - (R+1)(\ln R+\ln\ln R+1)
   \\ &=&R\ln\ln R-O(R)>0.
   \end{eqnarray*}  
   Otherwise we have e.g.\ $R+1\le c/(2\ln (c+1))$ and thus $c-(R+1)\ln(c+1)\ge c/2>0$, as desired. 
    On the other hand, for $c=R\ln R+R\ln\ln R$, we have very crudely that, say, $c+1\le \me R\ln R$ and thus
    $$
    c-(R+1)\ln(c+1) \le R\ln R+R\ln\ln R - (R+1)(\ln R+\ln\ln R+1)=(-1+o(1))R<0.
    $$
    Thus the largest root $c_0$ of~$\me^{c}=(c+1)^{R+1}$ indeed falls in the interval specified in~\eqref{eq:c0}.

  Thus
  \begin{eqnarray*}
  \mu &=& (c_0+1)\,\left(1+\frac{1}{c_0}\right)^{R}\ \le
  \ (R\ln R+2R\ln\ln R+1)\, \left(1+\frac{1}{R\ln R+R\ln\ln R}\right)^{R}
  \\ &\le& (R\ln R+2R\ln\ln R+1)\ \left(1+\frac{2}{\ln R}\right)\ <\ R\ln R+3R\ln\ln R.
  \end{eqnarray*}

 Now,  Corollary~\ref{cr:General} follows from Theorem~\ref{th:New}.\epf

\hide{
Given a concrete small $R$, it is easy to find numerical approximations to the smallest $\mu$ given by Theorem~\ref{th:New}. For example, the values of $\mu$ for $R=1,\dots,5$ rounded up in the third digit after the decimal point are $4.911$ (the same as in Theorem~\ref{th:main}), $9.268$, $14.036$, $19.113$ and $24.434$.
}
\hide{
\begin{verbatim}
r = 5; c =.; NSolve[E^c == (1 + c)^{r + 1}, c, Reals]

(c + 1)^(r + 1)/c^r /. %[[3]]
\end{verbatim}
}

\section{Concluding remarks}

One can re-write the proof of Lemma~\ref{lm:General} to also contain the conclusion of Lemma~\ref{lm:main}, with some minor changes for correctly handling the cases when $k<R$. However, the author feels that having a separate proof for the case $s=r+1$ is a good way to introduce the main ideas.

One can make the new lower bounds constructive, that is, for any fixed $R$, there is an algorithm that on input $(n,r)$ outputs a Tur\'an $(n,r+R,r)$-system in time polynomial in $n^r$. One has to replace a random subset $S\subseteq {[n]\choose k-R}$ by one constructed by the standard “conditioning method” (see e.g.~\cite{MotwaniRaghavan95ra}), in a very similar way as described in~\cite[Section IV]{KrivelevichSudakovVu03}.

Now that we know that $t(r+1,r)=O(1/r)$, the most intriguing remaining open question is whether  $t(r+1,r)=(1+o(1))/r$ as $r\to\infty$ or not.

\section*{Acknowledgements}

The author is very grateful to Oleg Verbitsky and Maksim Zhukovskii for drawing his attention to the paper~\cite{LenzRashtchianSiegelYaakobi21} and sharing their preliminary results exploiting analogies between Tur\'an systems and insertion covering codes. Also, the author would like to thank Xizhi Liu and the anonymous referees for their useful comments.
\medskip

For the purpose of open access, the author has applied a Creative Commons Attribution (CC-BY) licence to any Author Accepted Manuscript version arising from this submission.


\end{document}